\documentclass[12pt]{article}

\usepackage{latexsym}
\usepackage{amsmath}
\usepackage{amsfonts}
\usepackage{amssymb}
\usepackage{graphics}

\newcommand{\mb}{\mathbb}

\newcommand{\ot}{\otimes}

\newcommand{\wh}{\widehat}

\newcommand{\bean}{\begin{eqnarray}}
\newcommand{\eean}{\end{eqnarray}}
\newcommand{\bea}{\begin{eqnarray*}}
\newcommand{\eea}{\end{eqnarray*}}
\newcommand{\bsa}{\begin{subarray}{c}}
\newcommand{\esa}{\end{subarray}}
\newcommand{\bi}{\begin{itemize}}
\newcommand{\ei}{\end{itemize}}

\newtheorem{lemma}{Lemma}[section]
\newtheorem{thm}[lemma]{Theorem}

\def\ssp{\def\baselinestretch{1.0}\large\normalsize}

\title{ \bf Generalized Twisted Quantum Doubles and the McKay Correspondence}
\author{Chris Goff \\
University of the Pacific \\
Geoffrey Mason\thanks{Supported by NSA and NSF}\\
University of California at Santa Cruz \\
\\
Dedicated to Susan Montgomery}
\date{}
\begin{document}
\ssp
\maketitle

\abstract{We consider a class of quasi-Hopf algebras which we call \emph{generalized twisted quantum doubles}. They are abelian extensions $H = \mb{C}[\bar{G}] \bowtie \mb{C}[G]$ ($G$ is a finite group and $\bar{G}$ a homomorphic image), possibly twisted by a 3-cocycle, and are a natural generalization of the twisted quantum double construction of Dijkgraaf, Pasquier and Roche. We show that if $G$ is a subgroup 
of $SU_2(\mb{C})$ then $H$ exhibits an orbifold McKay Correspondence: certain fusion rules of $H$ define a graph with connected components indexed by conjugacy classes of $\bar{G}$, each connected 
component being an extended affine Diagram of type ADE whose McKay correspondent is 
the subgroup of $G$ stabilizing an element in the conjugacy class. This reduces to the original McKay Correspondence when $\bar{G} = 1$.}

\section{Introduction}
In an influential paper \cite{DPR},  Dijkgraaf, Pasquier and Roche introduced the twisted quantum double $D^{\omega}(G)$ of a finite group $G$. This is a quasi-Hopf algebra obtained by twisting the Drinfeld double $D(G)$ by a $3$-cocycle of $G$. They also suggested that the irreducible representations of $D^{\omega}(G)$ naturally correspond to the irreducible representations of a certain kind of conformal field theory (i.e. vertex operator algebra), namely a holomorphic orbifold\footnote{Here, $V$ is a holomorphic vertex operator algebra} $V^G$. Although this idea remains unproven at the level of mathematical rigor, it is almost certainly true. It is natural to ask if there are variants of the twisted double construction which might similarly correspond to other rational conformal field theories. Indeed, Dijkgraaf, Pasquier and Roche explicitly raised this question (loc. cit.) for the case of  theories with central charge $c=1$, where it is expected that most examples  arise from $G$-orbifolds where $G$ is a finite subgroup of $SU_2(\mathbb{C})$.

\bigskip
The purpose of the present  paper is to consider a class of quasi-Hopf algebras which we call
\emph{generalized twisted quantum doubles}.
They correspond to abelian extensions $H=\mathbb{C}[\bar{G}] \bowtie G$ where $\bar{G} = G/N$ with $N \unlhd G$, and may be twisted by a $3$-cocycle of $G$ that is inflated from a $3$-cocycle of 
$\bar{G}$. This produces a quasi-Hopf algebra that reduces to the DPR construction when $N=1$, and to the group algebra $\mathbb{C}[G]$ if $N=G$.  If $G$ is a finite subgroup of 
$SU_2(\mathbb{C})$ and $N$ has order at most $2$, we will see (Theorem \ref{McKayCorr}) that these quasi-Hopf algebras possess an 
\emph{orbifold McKay Correspondence}. By this we mean that we associate a graph $\mathcal{G}$ to $H$ whose vertices are the irreducible modules over $H$ and with edges given by the fusion rules for $H$ defined by tensoring with the canonical $2$-dimensional module for 
$\mathbb{C}[G]$ (which is also a module for $H$). Then the connected components of 
$\mathcal{G}$ are indexed by the conjugacy classes of $\bar{G}$; the connected component 
determined by $\bar{g} \in \bar{G}$ is an extended affine diagram of type ADE, namely the McKay correspondent of $C_G(\bar{g})$. Thus the graph $\mathcal{G}$ is indeed an orbifold version of the original correspondence of McKay \cite{Mc}, to which the construction reduces when $N=G$.

\bigskip
 We discuss several aspects of the orbifold McKay correspondence.\\
 1.  It is well-known in the physics literature that there is an ADE classification of $c=1$ rational conformal field theories (cf. \cite{MS} for more detail and further references). This arises in a manner that is rather different than our orbifold correspondence, but suggests  nevertheless that generalized twisted quantum doubles may indeed be related to certain orbifold conformal field theories.\\
2. The calculation of fusion rules of vertex operator algebras is usually difficult (cf.
 \cite{ADL} for an example of relevance to the present paper). The possibility of a McKay correspondence for fusion rules of vertex operator algebras suggests that, in some cases at least, there may be a more enlightening way to carry out the calculations.  The case when $N=1$, for example, shows (conjecturally) that \emph{all} holomorphic orbifolds 
 $V^G \ (G \subseteq SU_2(\mathbb{C}))$ exhibit an orbifold McKay correspondence. \\
 3. In current attempts to understand more general McKay Correspondences (e.g. \cite{IR}), the duality that exists between conjugacy classes and representations plays a r\^{o}le as an analog of  the duality between homology and cohomology. On the other hand, the graph $\mathcal{G}$ has a built-in duality: connected components are indexed by conjugacy classes, and
nodes of a connected component are indexed by irreducible modules over the stabilizer. \\

\bigskip
The paper is organized as follows. In Section $2$ we introduce our generalized twisted quantum doubles. We briefly develop some of their basic properties and the relationship between them and the usual twisted quantum doubles. In Section 3 we calculate fusion rules for generalized twisted quantum doubles. Some cases were already considered in \cite{DPR}. In Section 4 we describe the orbifold McKay correspondence. We will not treat the application of these results to orbifold conformal field theory here, but hope to return to this topic in the future.

\bigskip
Finally, it is a pleasure to record our debt to Susan Montgomery. She has had an enormous effect on our work as a friend, colleague, mentor, collaborator, and as the author of the most accessible book on Hopf algebras.

\section{Generalized twisted quantum doubles}
For background on Hopf algebras and related topics that we use here, the reader is referred to
\cite{Mo}.

\bigskip
Fix the following data and notation: $G$ is a finite group, $N \unlhd G$ a normal subgroup,
and $\bar{G} = G/N$. We
use the `bar convention' for elements in $\bar{G}$, i.e., if $g \in G$ then $\bar{g} = gN \in \bar{G}$.
$G$ acts by (right) conjugation on $\bar{G}$, so that $\bar{g}^x = \bar{g^x}= \bar{x}^{-1}\bar{g}\bar{x}.$ $\mathbb{C}[G]$ is the (complex) group algebra of $G$ and $\mathbb{C}[G]^*$ the dual of the group algebra. Let $\omega \in Z^3(\bar{G}, \mathbb{C}^*)$ be a multiplicative, normalized $3$-cocycle on $\bar{G}$ with
$\omega' =$ Infl$_{\bar{G}}^G\ \omega$ the inflation of $\omega$ to $G$. Thus
$\omega'(g, x, y) = \omega(\bar{g}, \bar{x}, \bar{y})$ for $g, x, y \in G$.  The associated $2$-cochains
$\theta$ and $\gamma$ are
\begin{eqnarray*}
\theta_{\bar{g}}(\bar{x}, \bar{y}) &=& \frac{\omega(\bar{g}, \bar{x}, \bar{y})\omega(\bar{x}, \bar{y}, \bar{g}^{xy})}{\omega(\bar{x}, \bar{g}, \bar{y}^g)},\\
\gamma_{\bar{g}}(\bar{x}, \bar{y}) &=& \frac{\omega(\bar{x}, \bar{y}, \bar{g})\omega(\bar{g}, \bar{x}^g, \bar{y}^g)}{\omega(\bar{x}, \bar{g}, \bar{y}^g)}.
\end{eqnarray*}

\medskip
For clarity we sometimes use the notation $\theta'$ and $\gamma'$ for the corresponding $2$-cochains associated with $\omega'$, so that $\theta'_g(x, y) = \theta_{\bar{g}}(\bar{x}, \bar{y})$
and $\gamma'_g(x, y) = \gamma_{\bar{g}}(\bar{x}, \bar{y})$.

\medskip
Define
\begin{eqnarray*}
D^{\omega}(G,  N) = \mathbb{C}[\bar{G}]^* \bowtie \mathbb{C}[G],
\end{eqnarray*}
where we use $\bowtie$ in place of $\otimes$ for notational convenience.
The product, coproduct, associator, counit, antipode and $\alpha$ and $\beta$ elements are defined as follows:

\begin{eqnarray*}
e(\bar{g}) \bowtie x. e(\bar{h}) \bowtie y &=& \delta_{\bar{g}^x, \bar{h}} \theta_{\bar{g}}(\bar{x}, \bar{y}) 
e(\bar{g}) \bowtie xy, \\
\Delta \ e(\bar{g})  \bowtie x &=& \sum_{\bar{a}\bar{b} = \bar{g}} \gamma_{\bar{x}} (\bar{a}, \bar{b})
e(\bar{a}) \bowtie x \otimes e(\bar{b}) \bowtie x, \\
\Phi &=& \sum_{\bar{g}, \bar{h}, \bar{k}} \omega(\bar{g}, \bar{h}, \bar{k})^{-1}
e(\bar{g}) \bowtie 1 \otimes e(\bar{h}) \bowtie 1 \otimes e(\bar{k}) \bowtie 1, \\
\epsilon \ e(\bar{g} \bowtie x) &=& \delta_{\bar{g}, 1}, \\
S \ e(\bar{g} \bowtie x) &=&  \theta_{\bar{g}^{-1}}( \bar{x}, \bar{x}^{-1}) \gamma_{\bar{x}}( \bar{g}, \bar{g}^{-1})^{-1} e(\bar{g}^{-\bar{x}}) \bowtie x^{-1}\\
\alpha &=& \mbox{Id} = \sum_{\bar{g}} \bar{g} \bowtie 1.\\
\beta &=& \sum_{\bar{g} \in  \bar{G}} \omega(\bar{g}, \bar{g}^{-1}, \bar{g}) \bowtie 1.
\end{eqnarray*}

This definition is, of course, modeled after the original twisted quantum double of Dijkgraaf-Pasquier-Roche (\cite{DPR}, \cite{Ka}), and the proof that it turns $D^{\omega}(G, N)$ into a quasiHopf algebra is the same. One simply has to make sure that the $\theta$- and $\gamma$-coefficients behave properly, and this is taken care of because the cocycle $\omega'$ on $G$ is inflated from $\omega$.
Note that $D^{\omega}(G, N)$ is also a cocentral abelian extension of Hopf algebras (cf. \cite{KMM}, Section 2). We call $D^{\omega}(G, N)$ a generalized (twisted) quantum double.
 If $N = G$ or $N = 1,$ then $D^{\omega}(G,  N)$ is
the group algebra $\mathbb{C}[G]$ or  the twisted quantum double $D^{\omega}(G)$
respectively.  

\bigskip
There are maps 
\begin{eqnarray}\label{qHdiagram}
D^{\omega}(G, N) \stackrel{\varphi}{\longrightarrow}  D^{\omega'}(G) \stackrel{\psi}{\longrightarrow} D^{\omega}(\bar{G}).
\end{eqnarray}
defined by
\begin{eqnarray}\label{varphidef}
&&\varphi:  e(\bar{g})\bowtie x  \mapsto \sum_{n \in N} e(gn) \bowtie x, \\
&&\psi: e(g) \bowtie x \mapsto e(\bar{g}) \bowtie \bar{x}.
\end{eqnarray}
Because $\omega'$ is inflated from $\omega$, it is evident that $\psi$ is a morphism of quasiHopf algebras. We assert that $\varphi$ is also a morphism of quasiHopf algebras.  We have
\begin{eqnarray*}
\varphi(e(\bar{g})\bowtie x). \varphi(e(\bar{h}) \bowtie y) &=&
\sum_{m, n \in N} e(gm) \bowtie x . e(hn) \bowtie y \\
&=&  \sum_{m, n \in N} \delta((gm)^x, hn) \theta'_{gm}(x, y) e(gm) \bowtie xy  \\
&=&  \delta(\bar{g}^x, \bar{h}) \theta_{\bar{g}}(\bar{x}, \bar{y})  \sum_{m \in N}     
e(gm) \bowtie xy \\
&=& 
 \delta(\bar{g}^x, \bar{h}) \theta_{\bar{g}}(\bar{x}, \bar{y}) 
\varphi(e(\bar{g}) \bowtie xy) \\
&=&\varphi(e(\bar{g})\bowtie x. e(\bar{h}) \bowtie y).
\end{eqnarray*}
So $\varphi$ preserves multiplication. Similarly,
\begin{eqnarray*}
\Delta \varphi(e(\bar{g})\bowtie x) &=&  \sum_{n \in N} \Delta e(gn) \bowtie x \\
&=& \sum_{n \in N} \sum_{ ab = gn} \gamma'_{x} ( a,  b)
e( a) \bowtie x \otimes e(b) \bowtie x \\
&=&  \sum_{ ab \in gN} \gamma_{\bar{x}} ( \bar{a},  \bar{b})
e( a) \bowtie x \otimes e(b) \bowtie x, 
\end{eqnarray*}
and
\begin{eqnarray*}
  \varphi(\Delta e(\bar{g})\bowtie x) &=&  \sum_{\bar{a}\bar{b} = \bar{g}} \gamma_{\bar{x}} (\bar{a}, \bar{b}) ( \varphi \otimes \varphi) (e(\bar{a}) \bowtie x \otimes e(\bar{b}) \bowtie x) \\
  &=&  \sum_{\bar{a}\bar{b} = \bar{g}} \gamma_{\bar{x}} (\bar{a}, \bar{b})
  \sum_{m, n \in N}  e(am) \bowtie x \otimes e(bn) \bowtie x \\
  &=& \sum_{a, b \in G, \bar{a}\bar{b} = \bar{g}} \gamma_{\bar{x}} (\bar{a}, \bar{b})
  e(a) \bowtie x \otimes e(b) \bowtie x,
  \end{eqnarray*}
so that $\varphi$ preserves comultiplication. We also check that $\varphi$ preserves
counits, associator, antipode, $\alpha$ and $\beta$ elements. Hence, 
$\varphi$ is indeed a morphism of quasiHopfalgebras. Let $H = im\ \varphi$.

\bigskip
 Now consider the left adjoint action 
  \begin{eqnarray*}
\mbox{ad}_l u (v) &=& \sum u_1 v  (Su_2),
\end{eqnarray*}
where we are using Sweeder notation $\Delta u =\sum u_1 \otimes u_2$.
Taking $u = e(h) \bowtie y, v = \varphi(e(\bar g) \bowtie x)$,
  \begin{eqnarray}\label{adlcalc}
\mbox{ad}_l u (v) &=& \sum_{a, b \in G, ab=h} \sum_{m \in N} 
\gamma'_y(a, b) \theta'_{b^{-1}}(y, y^{-1})^{-1}
e(a) \bowtie y. e(gm) \bowtie x.
  e(b^{-y}) \bowtie y^{-1} \notag \\
 &=&  \sum_{a, b \in G, ab=h} \sum_{m \in N}  \gamma'_y(a, b) \theta'_{b^{-1}}(y, y^{-1})^{-1}
  \theta'_a(y, x) \theta'_a(yx, y^{-1})
   \delta_{a^y, gm}  \delta_{a^{yx}, b^{-y}}
e(a) \bowtie x^{y^{-1}} \notag \\
&=&  \sum_{a, b \in G, ab=h} \sum_{m \in N}  \gamma'_y(a, b) \theta'_{b^{-1}}(y, y^{-1})^{-1}
  \theta'_a(y, x) \theta'_a(yx, y^{-1})
   \delta_{a^y, gm}  \delta_{a^{yx}, b^{-y}}
e(a) \bowtie x^{y^{-1}}.
   \end{eqnarray}
 
 Suppose that $H$ is a \emph{normal} subquasiHopf algebra of $D^{\omega'}(G)$. Then
 (\ref{adlcalc}) must lie in $im\ \varphi$ for all choices of $g, h, x$ and $y$. A summand can only be nonzero in case
 $a = (gm)^{y^{-1}}, b= a^{-yxy^{-1}}$ and $h = ab = ((gm)(gm)^{-x})^{y^{-1}}.$ If we therefore \emph{choose} $h = (gg^{-x})^{y^{-1}}$ (corresponding to $m=1$),  then the coefficient of $e(g^{y^{-1}}) \bowtie x^{y^{-1}}$ is a product of theta- and gamma-values, and in particular is nonzero. Because (\ref{adlcalc}) lies in $im\ \varphi$ then  the coefficients of $e(g^{y^{-1}}m) \bowtie x^{y^{-1}}$ are nonzero for each $m \in N$. (Indeed, all such coefficients are equal to that of $e(g^{y^{-1}} \bowtie x^{y^{-1}}$).) Therefore, the previous discussion shows that  $gg^{-x} = (gm)(gm)^{-x}$ is \emph{independent} of $m \in N$, and from this we readily find that $N \subseteq Z(G)$, the center of $G$.
 
 \bigskip
 Conversely, assume that $N \subseteq Z(G)$ with $h = (gg^{-x})^{y^{-1}}$. Setting $t=g^{y^{-1}}, u = x^{y^{-1}}$ and remembering that $\omega'$ is inflated from $\omega$, (\ref{adlcalc}) reads
   \begin{eqnarray}\label{adlcalc1}
\mbox{ad}_l u (v) 
&=&   \sum_{m \in N}  \gamma_{\bar{y}}(\bar{t}, \bar{t}^{-u}) 
\theta_{\bar{t}^u}(y, y^{-1})^{-1}
  \theta_{\bar{t}}(y, x) \theta_{\bar{t}}(yx, y^{-1}) e(tm) \bowtie u. 
   \end{eqnarray}
In particular, the coefficient of $e(tm) \bowtie u$ is independent of $m$, so that (\ref{adlcalc1}) indeed lies in $H$.
One similarly checks that the right adjoint ad$_ru(v)$ also lies in $H$. So we have proved
(cf. \cite{Mo}, P.33) that $H$ is a normal subquasiHopf algebra of $D^{\omega'}(G)$ if, and only if, $N \subseteq Z(G)$.

\bigskip
Let us continue to assume that $N \subseteq Z(G)$, and set 
$H^+ = H  \cap \mbox{ker}\ \epsilon, D = D^{\omega'}(G)$. 
Then $DH^+$ is a quasiHopf ideal in $D, D/DH^+$ a quasiHopf algebra, and the canonical projection
$D \rightarrow D/DH^+$ is a morphism of quasiHopf algebras. As a basis of
$D/DH^+$ we may take the (images of) either the  elements
$e(m) \bowtie 1$, or $e(1) \bowtie m$ for $m \in N$. Because $\omega$ is normalized and $\omega'$ is inflated from $\omega$, we have $e(m) \bowtie 1. e(n) \bowtie 1 = \delta_{m, n}e(m)\bowtie 1$. We thus obtain the following result.
\begin{lemma}\label{lemmaHnormal}
 Let $\varphi$ be as in (\ref{qHdiagram}), (\ref{varphidef}) and set 
 $D = D^{\omega'}(G), H = im\ \varphi$. Then 
 $\varphi$ is a morphism of quasiHopf algebras.
$H \subseteq D$ is a normal subquasiHopf algebra if, and only if,
$N \subseteq Z(G)$. In this case, there is an isomorphism of quasiHopf algebras $D/DH^+ \cong \mathbb{C}[N].$ $\hfill \Box$
\end{lemma}

\bigskip
 As is well-known (cf.\cite{Mo}, Chapter 7), the group algebra $\mathbb{C}[G]$ can be expressed as a 
\emph{crossed product} $\mathbb{C}[N] \#_{\sigma}\mathbb{C}[\bar{G}]$. Assuming for simplicity that
$N \subseteq Z(G)$ (the case of most interest to us), multiplication in the crossed product is
\begin{eqnarray*}
(m \# \bar{g})(n \#\bar{h}) = mn\sigma(\bar{g}, \bar{h}) \#\overline{gh}
\end{eqnarray*}
where $\sigma \in Z^2(\bar{G}, N)$ is a $2$-cocycle determined by the central extension
$1 \rightarrow N \rightarrow G \rightarrow \bar{G} \rightarrow 1$. Much as in the earlier calculations leading to Lemma \ref{lemmaHnormal}, we can check that there is an analogous description of
$D^{\omega}(G, N)$ as a crossed product using the \emph{same} $2$-cocycle. Precisely, we have
\begin{lemma}\label{lemmacrossprod} There is an isomorphism of quasiHopf algebras
\begin{eqnarray*}
D^{\omega}(G, N) \cong \mathbb{C}[N] \#_{\sigma} D^{\omega}(\bar{G}).
\end{eqnarray*}
$\hfill \Box$
\end{lemma}

By Lemma \ref{lemmacrossprod}, $D^{\omega}(G, N)$ is a central extension of the twisted quantum
double
$D^{\omega}(\bar{G})$, so that the representations of $D^{\omega}(G, N)$ are \emph{projective} representations of $D^{\omega}(\bar{G})$. From Lemma \ref{lemmaHnormal}, we can also use Clifford theory and obtain representations of $D^{\omega}(G, N)$ by restriction of representations of $D^{\omega'}(G)$.

\section{Fusion rules}
\bigskip
The simple modules over  $D^{\omega}(G, N)$ can be described as follows \cite{KMM}.
If $\{ \bar{g} \}$ is a set of representatives for the conjugacy classes of $\bar{G}$, 
the simple modules are 
\begin{eqnarray*}
\{ \mbox{Ind}_{C_G(\bar{g})}^G V \ | \ V \ \mbox{is a simple $\mathbb{C}^{\theta'_g}[C_G(\bar{g})]$-module} \}.
\end{eqnarray*}

If $V$ is a module over  $\mathbb{C}^{\theta'_g}[C_G(\bar{g})]$, we let $\widehat{V} = $Ind$^G_{C_G(\bar{g})}(V)$ denote the corresponding $D^{\omega}(G, N)$-module. Fix $\bar{g}$ and let $\chi_V$ be the character afforded by $V$.  The corresponding character $\wh{\chi}$ of $\wh{V}$ is then given by
\begin{eqnarray}\label{charval}
\wh{\chi}_{\wh{V}}(e(\bar{h}) \bowtie x) = \delta_{\bar{h}}^{y\bar{g}y^{-1}} \delta_{y^{-1}xy \in C_G(\bar{g})}  \theta'_{ygy^{-1}}(x,y) \theta_g^{'-1}(y,y^{-1}xy) \chi_V(y^{-1}xy).
\end{eqnarray}
The first $\delta$-function arises because the character value is zero unless $\bar{h}$ and $\bar{g}$ are conjugates.  We let $y$ be an element of $G$ which conjugates $\bar{g}$ to $\bar{h}$ in $\bar{G}$.  

\bigskip
In what follows, we let $U$ be a $C^{\theta'_f}_G(\bar{f})$-module, $V$ a $C^{\theta'_g}_G(\bar{g})$-module, and $W$ be a $C^{\theta'_h}_G(\bar{h})$-module where $\bar{f}, \bar{g}, \bar{h} \in \bar{G}$ are three of the conjugacy class representatives.
Let $\bar{J}$, $\bar{K}$, and $\bar{L}$ be the conjugacy classes of $\bar{G}$ that contain $\bar{f}$, $\bar{g}$, and $\bar{h}$, respectively. 
Define an inner product on the characters of $D^{\omega}(G, N)$ as follows:
\begin{equation} \label{innerprod}
\left< \wh{\chi}_{\wh{U}},\wh{\chi}_{\wh{V}}  \right> := \frac{1}{|G|} \sum_{\bar{k} \in \bar{G}} \sum_{x \in G} \delta_{\bar{k}}^{w\bar{f}w^{-1}} \delta_{\bar{k}}^{y\bar{g}y^{-1}} \delta_{x \in Q} \chi_U(w^{-1}xw)\overline{\chi_V(y^{-1}xy)}, 
\end{equation}
where $Q = wC_G(\bar{f})w^{-1} \cap yC_G(\bar{g})y^{-1}$.  Notice that the inner product is zero unless $\bar{f} = \bar{g}$.  In this case, 
\bea
\left< \wh{\chi}_{\wh{U}},\wh{\chi}_{\wh{V}}  \right> & =  & \frac{1}{|G|} \sum_{\bar{k} \in \bar{G}} \sum_{x \in G} \delta_{\bar{k}}^{y\bar{g}y^{-1}} \delta_{x \in yC_G(\bar{g})y^{-1}} \chi_U(y^{-1}xy) \overline{\chi_V(y^{-1}xy)} \\
& = & \frac{|\bar{J}|}{|G|} \sum_{y^{-1}xy \in C_G(\bar{g})} \chi_U(y^{-1}xy) \overline{\chi_V(y^{-1}xy)} \\
& = & \frac{1}{|C_G(\bar{g})|} \sum_{z \in C_G(\bar{g})} \chi_U(z) \overline{\chi_V(z)} \\
&=& \left< \chi_U, \chi_V \right>.
\eea
The last expression here is the usual inner product of characters of the group $C_G(\bar{g})$, although the characters may be projective. If the characters are \emph{ordinary}, the orthogonality relations for group characters implies that $\langle \chi_U, \chi_V \rangle = \delta_{U, V}$ in case $U, V$ are irreducible.  In general, we view the characters as ordinary characters of a \emph{covering group} of $C_G(\bar{g})$, and the orthogonality relations (applied to characters of the covering group) imply the same result.
We conclude  that  the irreducible characters of
$D^\omega(G, N)$ form an orthonormal basis with respect to (\ref{innerprod}).

\bigskip
Using the coproduct $\Delta$, we can derive the character of the tensor product module $\wh{V} \ot \wh{W}$.
\bea
\wh{\chi}_{\wh{V} \ot \wh{W}}(e(\bar{k})\bowtie x) & =  & \sum_{\bar{a}\bar{b} = \bar{k}} \gamma_{\bar{x}}(\bar{a},\bar{b}) \wh{\chi}_{\wh{V}}(e(\bar{a}) \bowtie x) \wh{\chi}_{\wh{W}}(e(\bar{b}) \bowtie x) \\
& = & \sum_{\bar{a}\bar{b} = \bar{k}}  \delta_{\bar{a}}^{y\bar{g}y^{-1}} \delta_{\bar{b}}^{z\bar{h}z^{-1}} \delta_{x \in C_G(\bar{a}) \cap C_G(\bar{b})}  \chi_V(y^{-1}xy) \chi_W(z^{-1}xz) \cdot \\
& \  & \cdot \  \gamma_{\bar{x}}(\bar{a},\bar{b}) \theta_{ygy^{-1}}(x,y) \theta_g^{-1}(y,y^{-1}xy)   \theta_{zhz^{-1}}(x,z) \theta_h^{-1}(z,z^{-1}xz).
\eea
 We obtain
\begin{eqnarray}\label{moreinnprod}
&&\left< \wh{\chi}_{\wh{V} \ot \wh{W}}, \wh{\chi}_{\wh{U}} \right> \notag \\
& = & \frac{1}{|G|} \sum_{\bar{k} \in \bar{G}} \sum_{x \in G} \sum_{\bsa \bar{a},\bar{b} \\ (\bar{a} \bar{b} = \bar{k}) \esa} \delta_{\bar{k}}^{w\bar{f}w^{-1}} \delta_{\bar{a}}^{y\bar{g}y^{-1}} \delta_{\bar{b}}^{z\bar{h}z^{-1}} \delta_{x \in C_G(\bar{a}) \cap C_G(\bar{b})} \chi_V(y^{-1}xy)\chi_W(z^{-1}xz) \overline{\chi_U(w^{-1}xw)} \cdot \notag \\
& \  &\ \ \ \ \ \ \ \ \ \ \ \ \ \ \ \ \ \ \ \ \  \cdot \  \gamma_{x}(a,b) \theta_{a}(x,y) \theta_g^{-1}(y,y^{-1}xy) \theta_{b}(x,z) \theta_h^{-1}(z,z^{-1}xz), 
\end{eqnarray}
that is 
\bea
\left< \wh{\chi}_{\wh{V} \ot \wh{W}}, \wh{\chi}_{\wh{U}} \right> & = & \frac{\delta_{\bar{J} \subseteq \bar{K}\bar{L}}}{|G|} \sum_{\bsa \bar{j} \in \bar{J} \\ (\bar{j} = w\bar{f}w^{-1}) \esa} \ \sum_{\bsa (\bar{a},\bar{b}) \in \bar{K} \times \bar{L} \\ (\bar{a}\bar{b} = \bar{j}) \\ (\bar{a} = y\bar{g}y^{-1}) \\ (\bar{b} = z\bar{h}z^{-1})  \esa} \ \sum_{x \in C_G(\bar{a}) \cap C_G(\bar{b})} \chi_V^{(y)}(x)\chi_W^{(z)}(x) \overline{\chi_U^{(w)}(x)} \cdot \\
& \  & \cdot \  \gamma_{x}(a,b) \theta_{a}(x,y) \theta_g^{-1}(y,y^{-1}xy) \theta_{b}(x,z) \theta_h^{-1}(z,z^{-1}xz),
\eea
where $\chi^{(t)}(x) = \chi(t^{-1}xt)$.

\bigskip
We specialize to the case where $\bar{h} = \bar{1}$.  Here, the previous displayed expression vanishes unless perhaps $\bar{f} = \bar{g}$. In this case,      
\bea
\left< \wh{\chi}_{\wh{V} \ot \wh{W}}, \wh{\chi}_{\wh{U}} \right> & = & \frac{1}{|G|} \sum_{\bsa \bar{k} \in \bar{J} \\(\bar{k} = w\bar{f}w^{-1}) \esa} \sum_{\bsa x \in C_G(\bar{k}) \\ (w^{-1}xw \in C_G(\bar{f})) \esa}   \chi_V(w^{-1}xw)\chi_W(x) \overline{\chi_U(w^{-1}xw)} \theta_{wfw^{-1}}(x,w) \theta_f^{-1}(w,w^{-1}xw). 
\eea
Since $W$ is a $G$-module, we can reindex the second summation over $w^{-1}xw$ to get
\begin{eqnarray}\label{innprod1}
&&\left< \wh{\chi}_{\wh{V} \ot \wh{W}}, \wh{\chi}_{\wh{U}} \right> \notag\\
& = & \frac{1}{|G|} \sum_{\bsa \bar{k} \in \bar{J} \\(\bar{k} = w\bar{f}w^{-1}) \esa} \sum_{t \in C_G(\bar{f})}   \chi_V(t)\chi_W(t) \overline{\chi_U(t)} \theta_{wfw^{-1}}(wtw^{-1},w) \theta_f^{-1}(w,t) \notag \\
& = &  \frac{1}{|G|} \sum_{t \in C_G(\bar{f})} \chi_V(t)\chi_W(t) \overline{\chi_U(t)} \left[ \sum_{w \in \{\textrm{coset reps of}\ C_G(\bar{f})\}}    \theta_{wfw^{-1}}(wtw^{-1},w) \theta_f^{-1}(w,t) \right].
\end{eqnarray}
We now have
\begin{lemma}\label{lemma3.1}
$\theta_{wfw^{-1}}(wtw^{-1},w) \theta_f^{-1}(w,t) =1$.
\end{lemma}
\begin{pf} Notice that the inner sum in (\ref{innprod1}) is independent of the modules
$U, V$ and $W$. Choose $W$ to be the trivial $1$-dimensional module for $G$, and let $U=V$ be the trivial $1$-dimensional module for $C_G(\bar{f})$. Then $\widehat{\chi}_W$ satisfies
$\widehat{\chi}_W(a) = \epsilon(a)$ for $a \in D^{\omega}(G, N)$ and  $\widehat{U} \otimes \widehat{W}
= \widehat{U}$. So in this case (\ref{innprod1}) reduces to
\begin{eqnarray}\label{innprod2}
1 =  \frac{1}{|G|} \sum_{t \in C_G(\bar{f})}  \left\{ \sum_{w \in \{\textrm{coset reps of}\ C_G(\bar{f})\}}    \theta_{wfw^{-1}}(wtw^{-1},w) \theta_f^{-1}(w,t) \right\}.
\end{eqnarray}
We may choose the values of $\omega$, and therefore also $\theta$, to be roots of unity. So the double sum in (\ref{innprod2}) is a sum of $|G|$ roots of unity equal to $|G|$. This means that each of the roots of unity is equal to $1$, thus proving the Lemma. $\hfill \Box$
\end{pf}

\bigskip
We can use Lemma \ref{lemma3.1} to simplify some of the earlier formulas. For example, the thetas may be removed in the character formula (\ref{charval}). Of more immediate concern is the fact that now reads
 (\ref{innprod1}) reads
\begin{eqnarray}\label{innprod3}
\left< \wh{\chi}_{\wh{V} \ot \wh{W}}, \wh{\chi}_{\wh{U}} \right> & = &  \frac{1}{|C_G(\bar{f})|} \sum_{t \in C_G(\bar{f})} \chi_V(t)\chi_W(t) \overline{\chi_U(t)} \notag \\
&=& \langle \chi_{V \otimes W}, \chi_U \rangle.
\end{eqnarray}
The last expression is the usual fusion rule for the twisted group algebra
$\mathbb{C}^{\theta'_f}[C_G(\bar{f})]$, and we regard $W$ as a module for $C_G(\bar{f})$ by restriction. We have therefore shown that fusion rules for $D^{\omega}(G, N)$ that involve a representation $W$ coming from $G = C_G(\bar{1})$ can be computed locally in the stabilizer $C_G(\bar{f})$.

\section{Orbifold McKay Correspondence}
From now on we specialize to the case that $G \subseteq SU_2(\mathbb{C})$. We further assume that
$G$ contains the (unique) involution $t = - I \in SU_2(\mathbb{C})$.
So we have the following commuting diagram with short exact rows.
     \begin{eqnarray*}
\begin{array}{ccccccccc}
 1& \longrightarrow&\mathbb{Z}_2& \longrightarrow&SU_2(\mathbb{C})&\longrightarrow & SO_3(\mathbb{R})&\longrightarrow& 1 \\
&&\parallel& & \uparrow& & \uparrow && \\
1&\longrightarrow& \mathbb{Z}_2& \longrightarrow& G & {\longrightarrow}&\bar{G}&
     \longrightarrow& 1 \\
\end{array}
\end{eqnarray*}
$G$ is a so-called \emph{binary polyhedral group}, defined abstractly via
\begin{eqnarray*}
\langle x, y, z \ | \ x^a = y^b = z^c = xyz \rangle,
\end{eqnarray*}
often denoted simply by $\langle a, b, c \rangle$. The maximal cyclic subgroups of $G$ have orders
$2a, 2b$ and $2c$. Apart from the cyclic case, 
$\langle a, b, c \rangle$ is finite only in the following cases: 
\begin{eqnarray*}
\langle a, b, c \rangle = \left \{ \begin{array}{l}
               \langle 2, 2, n \rangle \ (\mbox{$BD_{2n}$ = binary dihedral, order $4n$}),\\
                                                  \langle 2, 3, 3 \rangle  \ (\mbox{$SL_2(3)=$ binary tetrahedral, order $24$}),\\
                                             \langle 2, 3, 4 \rangle  \ (\mbox{$SL_2(3).2=$ binary octahedral, order $48$}),\\
                                             \langle 2, 3, 5 \rangle \ (\mbox{$SL_2(5)=$ binary icosahedral, order $120$}.
 \end{array}
\right.
\end{eqnarray*}

Each subgroup of $G$ is also isomorphic to one of these groups. 
The quotient groups $\bar{G}$ are cyclic of order $n$ ($Z_n$), dihedral of order $2n$ ($D_{2n}$),
$A_4, S_4$ and $A_5$ respectively. The $G$-orbits on $\bar{G}$ are just the conjugacy classes of $\bar{G}$.

\bigskip
It is  well-known (\cite{B}) that the second cohomology group (Schur multiplier) $H^2(G, \mathbb{C}^*)$ is trivial for each of these groups and subgroups. Thus, each of the $2$-cocycles
$\theta'_g$ is a coboundary, so that the twisted group algebra $\mathbb{C}^{\theta'_g}[C_G(\bar{g}]$ is isomorphic to the corresponding untwisted group algebra.  It follows that the number of irreducible modules for
$D^{\omega}(G, N)$ is independent of $\omega$ and is equal to the sum of the number of irreducibles for the stabilizers $C_G(\bar{g})$, one from each conjugacy class of $\bar{G}$. If  $\bar{g}$ has order greater than $2$ then $C_G(\bar{g})$ is a cyclic group of order $2b$ or $2c$. The individual cases are readily computed. Taking $N$ of order $2$, for example, there are the following possibilities:

\bigskip
\noindent
1. $G$ cyclic of order $2n$.
\begin{eqnarray*}
\begin{array}{c|c|c}
\bar{g} & C_G(\bar{g}) & \mbox{\# irreps of $C_G(\bar{g})$} \\ \hline
\mbox{any} & G & 2n
\end{array}
\end{eqnarray*}
Total $\#$ irreps $= 2n^2$.

\bigskip
\noindent
2. $G = \langle 2, 2, n \rangle$ ($n$ odd).
\begin{eqnarray*}
\begin{array}{c|c|c}
\bar{g} & C_G(\bar{g}) & \mbox{\# irreps of $C_G(\bar{g})$} \\ \hline
\bar{1} & G & n+3   \\
(n-1)/2 \ \mbox{classes} & Z_{2n} & 2n \\
\mbox{involution} & Z_4 & 4
\end{array}
\end{eqnarray*}
Total $\#$ irreps $= n^2+7$.

\bigskip
\noindent
3. $G = \langle 2, 2, n \rangle$ ($n$ even).
\begin{eqnarray*}
\begin{array}{c|c|c}
\bar{g} & C_G(\bar{g}) & \mbox{\# irreps of $C_G(\bar{g})$} \\ \hline
\bar{1} & G & n+3   \\
\mbox{central involution} & G & n+3 \\
(n-2)/2 \ \mbox{classes} & Z_{2n} & 2n \\
\mbox{involution ($2$ classes)} & BD_2 & 5
\end{array}
\end{eqnarray*}
Total $\#$ irreps $= n^2+16$.

\bigskip
\noindent
4. $G = \langle 2, 3, 3 \rangle$.
\begin{eqnarray*}
\begin{array}{c|c|c}
\bar{g} & C_G(\bar{g}) & \mbox{\# irreps of $C_G(\bar{g})$} \\ \hline
\bar{1} & G & 7   \\
 \mbox{(12)(34)} & BD_2 & 5 \\
\mbox{(123), (132)} & Z_6 & 6
\end{array}
\end{eqnarray*}
Total $\#$ irreps $= 24$.

\bigskip
\noindent
5. $G = \langle 2, 3, 4 \rangle.$ 
\begin{eqnarray*}
\begin{array}{c|c|c}
\bar{g} & C_G(\bar{g}) & \mbox{\# irreps of $C_G(\bar{g})$} \\ \hline
\bar{1} & G &  8  \\
\mbox{(12)} & BD_2 & 5 \\
\mbox{(12)(34)} & BD_4 & 7\\
\mbox{(1234)} & Z_8 & 8 \\
\mbox{(123)} & Z_6 & 6\\

\end{array}
\end{eqnarray*}
Total $\#$ irreps $= 34$.

\bigskip
\noindent
6. $G = \langle 2, 3, 5 \rangle.$
\begin{eqnarray*}
\begin{array}{c|c|c}
\bar{g} & C_G(\bar{g}) & \mbox{\# irreps of $C_G(\bar{g})$} \\ \hline
\bar{1} & G & 9   \\
(12)(34) & BD_2 & 5\\
(123) & Z_6 & 6 \\
(12345), (12354)& Z_{10} & 10 \\
\end{array}
\end{eqnarray*}
Total $\#$ irreps $= 40$.

\bigskip
Let $W$ be the canonical $2$-dimensional irreducible $G$-module affording
the embedding  $G \rightarrow SU_2(\mathbb{C})$ if $G$ is \emph{not} cyclic.
If $G$ is cyclic, we take $W$ to be the direct sum of a $1$-dimensional faithful $G$-module and its dual.
As a
$C_G(\bar{1})$-module, we may, and shall, consider $W$ as a module over $D^{\omega}(G, N)$. In the notation of Section 3, $W = \widehat{W}$. 

\bigskip
Following the original construction of McKay (\cite{Mc}), introduce a graph $\mathcal{G}$
whose vertices are the irreducible modules over $D^{\omega}(G, N)$. If $\widehat{U}, \widehat{V}$
are two vertices, we connect them by $ \langle \widehat{\chi}_{\widehat{U}\otimes\widehat{W} }, \widehat{\chi}_{\widehat{V}} \rangle$ edges.
We saw in the last Section that $ \langle \widehat{\chi}_{\widehat{U}\otimes\widehat{W}}, \widehat{\chi}_{\widehat{V}} \rangle = 0$ unless perhaps $U$ and $V$ are both irreducible modules for
some $\mathbb{C}^{\theta'_f}[C_G(\bar{f})]$, in which case $ \langle \widehat{\chi}_{\widehat{U}\otimes\widehat{W} }, \widehat{\chi}_{\widehat{V}} \rangle= \langle \chi_{U \otimes W}, \chi_V \rangle$.
In particular, vertices indexed by modules in \emph{distinct} stabilizers are \emph{not} connected.
On the other hand, if $\bar{f}$ is fixed and $U, V$ are irreducible modules for 
$\mathbb{C}^{\theta'_f}[C_G(\bar{f})]$,  the  multiplicities $ \langle \chi_{U \otimes W}, \chi_V \rangle$ are \emph{precisely} those which result from the McKay procedure (loc. cit.) applied to the stabilizer $C_G(\bar{f})$. Applying the McKay Correspondence in these cases, we arrive at our main result, which we state as follows.
\begin{thm}\label{McKayCorr} Let the notation be as above. The connected components of the graph $\mathcal{G}$ 
are indexed by the conjugacy classes of $\bar{G}$. The connected component associated to
the class of $\bar{g}\in \bar{G}$ is the affine Dynkin diagram $\Phi_{\bar{g}}$ which is the McKay correspondent of the stabilizer $C_G(\bar{g})$. If $\bar{g}$ as order greater than $2$ then
$\Phi_{\bar{g}}$ is of type $\hat{A}_{2j-1}$ for some $j$. $\hfill \Box$
\end{thm}

The graph $\mathcal{G}$ is described in each case by the data given earlier in this Section. Here is a picture of the case when $G = \langle2, 3, 4\rangle$ and $N$ has order $2$.

\begin{center}
\resizebox{8cm}{!}{\includegraphics{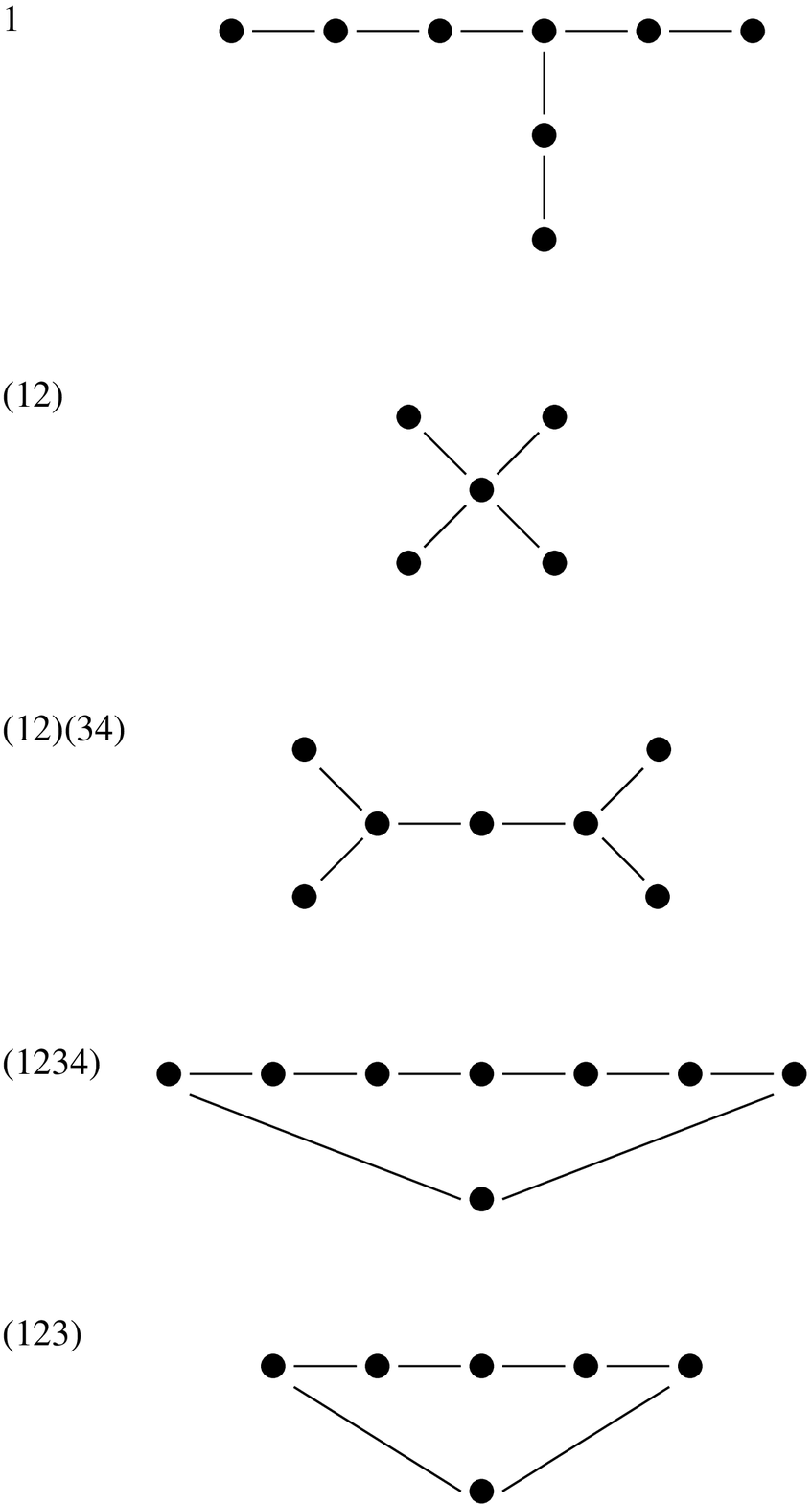}}
\end{center}

\end{document}